\documentclass{amsart}

\usepackage{amssymb,amscd,amsmath, latexsym}
\usepackage[all]{xy}
\usepackage{hyperref}
  \newtheorem{proposition}{Proposition}[section]

  \newtheorem{corollary}[proposition]{Corollary}
  \newtheorem{theorem}[proposition]{Theorem}

  \theoremstyle{definition}
  \newtheorem{definition}[proposition]{Definition}

  \theoremstyle{remark}
  \newtheorem{remark}[proposition]{Remark}

\author{Marino Gran and Diana Rodelo}

\address{Insitut de Math\'ematique Pure et Appliqu\'ees, Universit\'e Catholique de Louvain, B-1348
Louvain-la-Neuve, Belgium \\
E-mail address: marino.gran@uclouvain.be \\}
\address{
Departamento de Matem\'atica, Faculdade de Ci\^{e}ncias e Tecnologia, \\
Universidade do Algarve, Campus de Gambelas, 8005-139 Faro, Portugal \\
Centro de Matem\'atica da Universidade de Coimbra, 3001-454 Coimbra, Portugal\\
E-mail address: drodelo@ualg.pt}

\thanks{Research supported by the F.N.R.S. grant \emph{Cr\'edit aux chercheurs} 1.5.016.10F, and by FCT/Centro de Matem\'atica da Universidade de Coimbra}

\title{On the characterization of J\'onsson-Tarski  and of subtractive varieties}%



\begin{document}

\newdir{ >}{% end of arrow for the monos
  @{}*!/-10pt/@{>} }
\newdir{ |>}{% end of arrow for the kernels
  @{}*!/-5.5pt/@{|}*!/-10pt/:(1,-.2)@^{>}*!/-10pt/:(1,+.2)@_{>} }

\newcommand{\EE}{ \ensuremath{\mathbb {E}} }
\newcommand{\CC}{ \ensuremath{\mathbb {C}} }

\newcommand{\map}[2]{ \ensuremath{ \xymatrix@1@C=15pt{ #1 \ar[r] & #2 } } }
\newcommand{\mono}[2]{ \ensuremath{ \xymatrix@1@C=15pt{ #1 \ar@{ >->}[r] & #2 } } }
\newcommand{\regepi}[2]{ \ensuremath{ \xymatrix@1@C=15pt{ #1 \ar@{>>}[r] & #2 } } }
\newcommand{\Span}[5]{ \ensuremath{ \xymatrix@1@C=15pt{ #1 & #3 \ar[l]_-{#2} \ar[r]^-{#4} & #5 } } }
\newcommand{\punctualspan}[7]{ \ensuremath{ \xymatrix@1@C=15pt{ #1 \ar[r]_-{#3} &
#4 \ar@<-3pt>[l]_-{#2} \ar@<3pt>[r]^-{#5} & #7 \ar[l]^-{#6} } } }
\newcommand{\rightpunctualspan}[7]{ \ensuremath{ \xymatrix@1@C=15pt{ #1 &
#4 \ar@<-3pt>[l]_-{#2} \ar@<3pt>[r]^-{#5} & #7 \ar[l]^-{#6} } } }
\newcommand{\longspan}[5]{ \ensuremath{ \xymatrix@1@C=30pt{ #1 & #3 \ar[l]_-{#2} \ar[r]^-{#4} & #5 } } }
\newcommand{\longpunctualspan}[7]{ \ensuremath{ \xymatrix@1@C=30pt{ #1 \ar[r]_-{#3} &
#4 \ar@<-3pt>[l]_-{#2} \ar@<3pt>[r]^-{#5} & #7 \ar[l]^-{#6} } } }
\newcommand{\longrightpunctualspan}[7]{ \ensuremath{ \xymatrix@1@C=30pt{ #1 &
#4 \ar@<-3pt>[l]_-{#2} \ar@<3pt>[r]^-{#5} & #7 \ar[l]^-{#6} } } }
\newcommand{\prd}[2]{ \ensuremath % horizontal style morphism for products
	\langle #1,#2 \rangle}
\newcommand{\cop}[2]{\ensuremath % vertical style morphism for sums
	\left< \begin{smallmatrix} #1 \\ #2 \end{smallmatrix} \right>}
	%\fontsize{8}{10} \left \langle \!\begin{array}{c} #1 \vspace{-4pt}\\ #2 \end{array}  \!\right \rangle}
\newcommand{\add}{\ensuremath \textbf{+}}
	
\hyphenation{e-quiv-a-len-ces co-ker-nels gen-er-al-ised ex-act-ness ex-ten-sion Ja-ne-lid-ze pro-jec-tive co-ker-nel
group-oid group-oids sub-trac-tive w-sub-trac-tive}

\maketitle

\begin{center} \emph{Dedicated to Madame Andr\'ee Ehresmann on the occasion of her 75th birthday} \end{center}
\vspace{10pt}

%%%%%%%%%%%%%%%%%%%%%%%%%%%%%%%%%%%%%%%  ABSTRACT  %%%%%%%%%%%%%%%%%%%%%%%%%%%%%%%%%%%%%%%%%%%%%%%%%%%%%%%%%%%%%%%%%%%%
\begin {abstract}
We investigate unital, subtractive and strongly unital regular categories with enough projectives and give characterizations of
their projective covers. The categorical equation \emph{strongly unital = unital + subtractive} is explored: this leads to 
 proofs of their varietal characterizations in terms of the categorical properties of the corresponding algebraic theories.
\end {abstract}

%%%%%%%%%%%%%%%%%%%%%%%%%%%%%%%%%%%%%%%%%%%%%%%%%%%%%%%%%%%%%%%%%%%%%%%%%%%%%%%%%%%%%%%%%%%%%%%%%%%%%%%%%%%%%%%%%%%%%%%
%%%%%%%%%%%%%%%%%%%%%%%%%%%%%%%%%%%%%%%  INTRODUCTION  %%%%%%%%%%%%%%%%%%%%%%%%%%%%%%%%%%%%%%%%%%%%%%%%%%%%%%%%%%%%%%%%
%%%%%%%%%%%%%%%%%%%%%%%%%%%%%%%%%%%%%%%%%%%%%%%%%%%%%%%%%%%%%%%%%%%%%%%%%%%%%%%%%%%%%%%%%%%%%%%%%%%%%%%%%%%%%%%%%%%%%%%
\section*{Introduction} 

Over the past years, there have been several investigations devoted to the construction of ``free'' categories of a specific type, a
 classical instance of this being the free abelian category over a (pre)additive one \cite{RAC, Adelman}. 

Abelian categories are, in particular, exact categories \cite{EC}: this naturally led to the study of the free exact category over a category with finite limits, that was first described in \cite{FECLEO}. In particular, this result was useful to obtain a conceptual construction in two steps of the free abelian category mentioned above. It was later observed that, for the free exact completion, the assumption that the original category had finite \emph{limits} was never fully used, and only the existence of finite \emph{weak limits} was actually needed (recall that for a weak limit the uniqueness requirement in the definition of a limit is simply dropped). The exact completion over a category with weak finite limits \cite{REC}, called \emph{weakly lex} in the following, is general enough to cover many further examples, such as any algebraic category. 

An interesting feature of free exact categories over weakly lex ones is that they always possess enough regular projectives, a property which is quite useful in homological algebra.
One actually knows even more: an exact category is the exact completion of a weakly lex category if and only if it has enough regular projectives and, in this case, it is the completion of its full subcategory of regular projectives. More generally, any exact category with enough projectives $\EE$ is the exact completion of any of its \emph{projective covers}, that is of any full subcategory $\CC$ of $\EE$ whose objects are regular projectives and that has, moreover, the property that for any object $A$ in $\EE$ there is a regular epimorphism $\regepi{\tilde{A}}{A}$ with $\tilde{A}$ in $\CC$.
For example, any variety of universal algebras turns out to be the exact completion of its full subcategory of free algebras. Similar considerations also apply to the \emph{free regular completion} over a weakly lex category (see \cite{Vit} and references therein).
Many results on abelian categories can be recovered by using these completions, and the proofs are general and simple (see \cite{ECRAC}). 

An important aspect in the study of these completions concerns the possibility of characterizing those categories that appear as \emph{projective covers} of a free exact category. For instance, already in \cite{RAC} it was shown that a category is equivalent to a projective cover of an abelian category if and only if it is preadditive with weak finite products and weak kernels.  
Similar characterizations have been later established for the projective covers of categories that are Mal'tsev \cite{ECRAC}, protomodular, semi-abelian \cite{G}, extensive \cite{ECHC}, (locally) cartesian closed categories \cite{LCCEC, CCEC} and toposes \cite{Menni}.

The aim of this work is to extend the above list, by giving characterizations of the projective covers of regular categories $\EE$ which are unital, subtractive \cite{SC} and strongly unital \cite{MCFPO}. The investigation of these specific properties is inspired from the ``categorical equation''
\begin{center} \emph{strongly unital = unital + subtractive}. \end{center} 
When $\EE$ is a variety of universal algebras, our results isolate those properties of the categories of free algebras that ``force'' $\EE$ to be a unital, subtractive or strongly unital variety, respectively.
As a consequence we obtain some new proofs of their varietal characterizations in terms of the categorical properties of the corresponding algebraic theories (in the sense of \cite{Law}).
The ideas behind our results are closely related to those which appear in the joint work of D. Bourn and Z. Janelidze (see \cite{Zurab} and the references there), and it would be interesting to
investigate the precise connections with their work in the future.

%%%%%%%%%%%%%%%%%%%%%%%%%%%%%%%%%%%%%%%%%%%%%%%%%%%%%%%%%%%%%%%%%%%%%%%%%%%%%%%%%%%%%%%%%%%%%%%%%%%%%%%%%%%%%%%%%%%%%%%
%%%%%%%%%%%%%%%%%%%%%%%%%%%%%%%%%%%%%%%  Regular completions  %%%%%%%%%%%%%%%%%%%%%%%%%%%%%%%%%%%%%%%%%%%%%%%%%%%%%%%%%
%%%%%%%%%%%%%%%%%%%%%%%%%%%%%%%%%%%%%%%%%%%%%%%%%%%%%%%%%%%%%%%%%%%%%%%%%%%%%%%%%%%%%%%%%%%%%%%%%%%%%%%%%%%%%%%%%%%%%%%
\section{Projective covers} 
\label{Projective covers} 
When $\CC$ is a full subcategory of a category $\EE$ one says that $\CC$ is a \emph{projective cover} of $\EE$ if two conditions are satisfied: 1) any object in $\CC$ is regular projective in $\EE$; 2) for any object $A$ in $\EE$ there exists a $\CC$-cover $\tilde{A}$ of $A$, that is an object $\tilde{A}$ in $\CC$ and a regular epimorphism $\regepi{\tilde A}{A}$.

For example, if $\EE$ is a (quasi)variety of universal algebras and $\CC$ its full subcategory of free algebras, then $\CC$ is a projective cover of the (regular) exact category $\EE$.

Recall that a finitely complete category is called \emph{regular} \cite{EC} when it has a pullback-stable (regular epimorphism, monomorphism) factorization system.

%%%%%%%%%%%%%%%%%%%%%%%%%%%%%%%%%%%%%%%  REMARK: PROPOSITION 4 IN REC  %%%%%%%%%%%%%%%%%%%%%%%%%%%%%%%%%%%%%%%%%%%%%%%%
\begin{remark}
\label{Remark: Proposition 4 in REC}
Let $\CC$ be a projective cover of a category $\EE$. As a consequence of Proposition 4 in \cite{REC}, we have that
$\CC$ is \emph{weakly lex} (i.e. it
has weak finite limits) whenever $\EE$ is weakly lex. For instance, if $X,Y$ are objects in $\CC$, $\Span{X}{}{E}{}{Y}$ is a weak product in
$\EE$ and $\regepi{\tilde{E}}{E}$ a $\CC$-cover of $E$, then $\Span{X}{}{\tilde{E}}{}{Y}$ is a weak product in $\CC$. On the other hand, if
$\EE$ is a regular category and $\Span{X}{}{C}{}{Y}$ is a weak product in $\CC$, then the factorization
$\regepi{C}{X\times Y}$ is a regular epimorphism. Similar remarks apply also to all (weak) finite limits.
\end{remark}

Note that, if $\CC$ is a pointed projective cover of a category $\EE$ with kernel pairs, then $\EE$ is also pointed. Let us recall the argument: if $P$ is the zero object of $\CC$, then it must be the initial object of $\EE$, since, for any object $A$ of
$\EE$, we can consider the $\CC$-cover $a:\regepi{\tilde A}{A}$ and the unique morphism $i:\map{P}{\tilde A}$ in $\CC$;
this gives a morphism $\map{P}{A}$. For the uniqueness, suppose that there exists another morphism $x:\map{P}{A}$. The fact that $P$ is projective gives a morphism $j:\map{P}{\tilde A}$ such that $x=a\cdot j$, and $i=j$ since $P$ is an initial object in $\CC$. As for being a terminal object in $\EE$, first we notice that there exists a unique morphism
$t:\map{\tilde A}{P}$ in $\CC$. Next, we consider the $\CC$-cover $r:\map{\tilde R}{R[a]}$ of the kernel pair $(R[a], a_1, a_2)$ of $a$.
Then, we have the equality $t\cdot a_1\cdot r=t\cdot a_2\cdot r$ of two morphisms onto the terminal object of $\CC$. This
gives a (unique) morphism $\pi:\map{A}{P}$ such that $\pi\cdot a=t$ since $a$ is the coequalizer of its kernel pair in
$\EE$. The uniqueness of the morphism $\map{A}{P}$ comes from the fact that precomposition with the epimorphism $a$ is the
unique morphism $t$ in $\CC$.

In the next sections we shall give characterizations of the projective covers of regular categories which are unital, subtractive and strongly unital. 

%%%%%%%%%%%%%%%%%%%%%%%%%%%%%%%%%%%%%%%%%%%%%%%%%%%%%%%%%%%%%%%%%%%%%%%%%%%%%%%%%%%%%%%%%%%%%%%%%%%%%%%%%%%%%%%%%%%%%%%
\section{Unital categories}
\label{Unital categories}

We begin this section by recalling the main properties of unital categories \cite{MCFPO}. We then give a ``weak'' version of such
categories in order to characterize the projective covers of regular unital categories. In the presence of finite coproducts, we
shall obtain a simple proof of the characterization of algebraic (quasi)varieties which are unital.

In a category with pullbacks, strong and extremal epimorphisms coincide. If, moreover, the category is regular, then the notions of regular, strong and extremal epimorphisms coincide.

To simplify notation, we write $1$ (instead of $1_X$) to denote the identity morphism on any object $X$ and just $0$ to denote the null morphism from an object $X$ to an object $Y$, when the category is pointed.

%%%%%%%%%%%%%%%%%%%%%%%%%%%%%%%%%%%%%%%  DEFINITION: UNITAL CATEGORY  %%%%%%%%%%%%%%%%%%%%%%%%%%%%%%%%%%%%%%%%%%%%%%%%%
\begin{definition}\emph{(\cite{MCFPO})
\label{Unital} 
A category $\EE$ is called \emph{unital} when it is pointed, has finite limits and, for every pair of objects $X,Y$ in
$\EE$, the pair of morphisms $$\xymatrix@1@C=30pt{X \ar[r]^-{\langle 1,0 \rangle} & X\times Y &  Y \ar[l]_-{\langle 0,1 \rangle}}$$ is jointly strongly
epimorphic.}
\end{definition}

%%%%%%%%%%%%%%%%%%%%%%%%%%%%%%%%%%%%%%%  PUNCTUAL SPAN  %%%%%%%%%%%%%%%%%%%%%%%%%%%%%%%%%%%%%%%%%%%%%%%%%%%%%%%%%%%%%%%
Let $\EE$ be a pointed category. A span $\Span{X}{f}{Z}{g}{Y}$ in $\EE$ is called
\emph{punctual} if there exist morphisms $s:\map{X}{Z}$ and $t:\map{Y}{Z}$ such that $f\cdot s=1, g\cdot t=1, f\cdot
t=0$ and $g\cdot s=0$. For example, for any pair of objects $X,Y$ in $\EE$, we have a punctual span
$\longpunctualspan{X}{\pi_1}{\langle 1,0 \rangle}{X\times Y}{\pi_2}{\langle 0,1 \rangle}{Y}$. Unital categories can be characterized through such spans
as follows (see Theorem 1.2.12 of \cite{BB}):

%%%%%%%%%%%%%%%%%%%%%%%%%%%%%%%%%%%%%%%  THEOREM: CHARACTERIZATION OF UNITAL CATEGORIES  %%%%%%%%%%%%%%%%%%%%%%%%%%%%%%
\begin{theorem}
\label{Theorem: characterization of unital cats} A pointed category with finite limits $\EE$ is unital if and only
if for every punctual span $\punctualspan{X}{f}{s}{Z}{g}{t}{Y}$ in $\EE$, the induced factorisation $\langle f,g \rangle:\regepi{Z}{X\times
Y}$ is a strong epimorphism.
\end{theorem}

We introduce the following definition, which is suggested by the observations made in Remark~\ref{Remark: Proposition 4 in REC}.

%%%%%%%%%%%%%%%%%%%%%%%%%%%%%%%%%%%%%%%  DEFINITION: w-UNITAL CATEGORY  %%%%%%%%%%%%%%%%%%%%%%%%%%%%%%%%%%%%%%%%%%%%%%%
\begin{definition}
\label{Definition: w-unital cat}
\emph{A category $\CC$ is called \emph{w-unital} when it is pointed, has weak finite limits and for every punctual span
$\punctualspan{X}{f}{s}{Z}{g}{t}{Y}$, the span $\Span{X}{f}{Z}{g}{Y}$ is a weak product.}
\end{definition}

It is easy to see that any finitely complete w-unital category is necessarily unital, since the factorization to the
``true'' product $\prd{f}{g}:\regepi{Z}{X\times Y}$ is a split epimorphism, thus a strong epimorphism. Similar observations can be made also for subtractive and strongly unital categories, which shall be analysed in the next sections.

%%%%%%%%%%%%%%%%%%%%%%%%%%%%%%%%%%%%%%%  PROPOSITION: CC w-UNITAL <=> EE UNITAL  %%%%%%%%%%%%%%%%%%%%%%%%%%%%%%%%%%%%%%
\begin{proposition}
\label{Proposition: CC w-unital <=> EE unital}
Let $\CC$ be a pointed projective cover of a regular category $\EE$. Then $\EE$ is unital if and only if $\CC$ is
w-unital.
\end{proposition}
\proof Suppose that $\EE$ is unital. Since $\CC$ is a projective cover of a regular category $\EE$, which has finite
limits, then $\CC$ has weak finite limits (Remark~\ref{Remark: Proposition 4 in REC}). Now, let $\punctualspan{\tilde
X}{\tilde f}{\tilde s}{\tilde Z}{\tilde g}{\tilde t}{\tilde Y}$ be a punctual span in $\CC$. Since it is also a
punctual span in $\EE$, which is unital, then $\prd{\tilde f}{\tilde g}:\regepi{\tilde Z}{\tilde X\times \tilde Y}$ is a
regular epimorphism in $\EE$. In other words, it represents a $\CC$-cover of the true product, so it is a weak product in $\CC$ (Remark~\ref{Remark: Proposition 4 in REC}).

\noindent Conversely, $\EE$ is pointed, because $\CC$ is pointed (see Section~\ref{Projective covers}), and it has finite limits. Let
$\punctualspan{X}{f}{s}{Z}{g}{t}{Y}$ be a punctual span in $\EE$. We consider the $\CC$-covers $x:\regepi{\tilde X}{X}$
and $y:\regepi{\tilde Y}{Y}$, form the pullback $U$ of $x\times y$ and $\langle f,g \rangle$ and take its $\CC$-cover
$\alpha:\regepi{\tilde Z}{U}$
$$
  \xymatrix@1@R=25pt{\tilde Z \ar@{>>}[r]^-{\alpha} \ar[dr]_-{\langle\tilde f=f'\cdot \alpha,\, \tilde g=g'\cdot \alpha \rangle}
  & U \ar@{>>}[r]^-{u} \ar[d]|(.4){\langle f',g' \rangle\;\;} \ar@{}[dr]|(.25){\lrcorner} & Z \ar[d]^-{\langle f,g \rangle} \\
  & \tilde X\times \tilde Y \ar@{>>}[r]_-{x\times y} & X\times Y. }
$$
Note that, $\punctualspan{\tilde X}{f'}{s'}{U}{g'}{t'}{\tilde Y}$ is a punctual span in $\EE$, for $s'=(\prd{1}{0},s\cdot
x)$ and $t'=(\prd{0}{1},t\cdot y)$. Since $\tilde{X}$ and $\tilde{Y}$ are projective and $\alpha$ is a regular epimorphism,
then there exists morphisms $\tilde{s}$ and $\tilde{t}$ such that $s'=\alpha\cdot \tilde{s}$ and $t'=\alpha\cdot \tilde{t}$.
We obtain a punctual span $\punctualspan{\tilde X}{\tilde f}{\tilde s}{\tilde Z}{\tilde g}{\tilde t}{\tilde Y}$ in
$\CC$. Since $\CC$ is w-unital, then $(\tilde f, \tilde g)$ is a weak product and, consequently, $\prd{\tilde f}{\tilde
g}:\regepi{\tilde Z}{\tilde X\times \tilde Y}$ is a regular epimorphism (Remark~\ref{Remark: Proposition 4 in REC}).
Finally, the equality $\prd{f}{g}\cdot u\cdot \alpha=x\times y\cdot \prd{\tilde{f}}{\tilde{g}}$ allows us to conclude that
$\prd{f}{g}$ is a regular epimorphism, thus a strong epimorphism and, therefore, $\EE$ is unital.
$\Box$\vspace{5pt}

When working in a context that admits binary coproducts, (w-)unital categories can be characterized through a special
type of punctual span (see Proposition 1.2.18 in \cite{BB} for unital categories):

%%%%%%%%%%%%%%%%%%%%%%%%%%%%%%%%%%%%%%%  PROPOSITION: CHARACTERIZATION OF w-UNITAL CATEGORIES WITH +  %%%%%%%%%%%%%%%%%%%
\begin{proposition}
\label{Proposition: characterization of w-unital cats with +}
A pointed weakly lex category with binary coproducts $\CC$ is w-unital if and only, for every pair of objects $X,Y$ in
$\CC$, the canonical span $$\longspan{X}{\cop{1}{0}}{X+ Y}{\cop{0}{1}}{Y}$$ is a weak product.
\end{proposition}
\proof When $\CC$ is w-unital, the span  $\longspan{X}{\cop{1}{0}}{X+ Y}{\cop{0}{1}}{Y}$ is a weak product, since the span $\longpunctualspan{X}{\cop{1}{0}}{i_1}{X+Y}{\cop{0}{1}}{i_2}{Y}$ is punctual.

\noindent Conversely, consider a punctual span $\punctualspan{X}{f}{s}{Z}{g}{t}{Y}$ and arbitrary morphisms
$x:\map{A}{X}$, $y:\map{A}{Y}$. By assumption, there is a morphism $w:\map{A}{X+ Y}$ such that $\cop{1}{0}\cdot w=x$ and
$\cop{0}{1}\cdot w=y$. So, there exists a morphism $\cop{s}{t}\cdot w:\map{A}{Z}$ such that $f\cdot \cop{s}{t}\cdot w=x$ and $g\cdot
\cop{s}{t}\cdot w=y$, proving that $(f,g)$ is a weak product.
$\Box$\vspace{5pt}

%%%%%%%%%%%%%%%%%%%%%%%%%%%%%%%%%%%%%%%  JONSSON-TARSKI VARIETIES  %%%%%%%%%%%%%%%%%%%%%%%%%%%%%%%%%%%%%%%%%%%%%%%%%%%%
Recall that a variety of universal algebras is called a \emph{J\'onsson-Tarski} variety when its theory contains a unique constant 0 and a binary
operation + such that $x+0=x=0+x$ \cite{JT}. An \emph{(internal) unital coalgebra} in any category $\EE$ with binary coproducts is defined as an object $X$ equipped with a morphism $\add:\map{X}{X+X}$ such that $\cop{1}{0} \cdot \add = 1 = \cop{0}{1} \cdot \add$.

It is known that a variety is J\'onsson-Tarski if and only if it is a
unital category (Theorem 1.2.15 of \cite{BB}). This result is illustrated in the next theorem, where we give a new proof based on the structural property of the category of free algebras (given in Proposition \ref{Proposition: characterization of w-unital cats with +}).

Given a pointed algebraic variety $\EE$ and its full subcategory $\CC$ of free algebras, we are in the situation where
$\CC$ is a pointed projective cover with binary coproducts of the regular category $\EE$.

%%%%%%%%%%%%%%%%%%%%%%%%%%%%%%%%%%%%%%%  THEOREM: JONSSON-TARSKI AND (w-)UNITAL  %%%%%%%%%%%%%%%%%%%%%%%%%%%%%%%%%%%%%%
\begin{theorem}
\label{Theorem: Jonsson-Tarski and (w-)unital}
Let $\EE$ be an algebraic variety and $\CC$ its full subcategory of free algebras. Then the following statements are
equivalent:
\begin{itemize}
  \item[1.] $\EE$ is a unital category;
  \item[2.] $\CC$ is a w-unital category;
  \item[3.] the free algebra $X=F(1)$ on one generator is a unital coalgebra;
  \item[4.] $\EE$ is a J\'onsson-Tarski variety.
\end{itemize}
\end{theorem}
\proof  $1. \Leftrightarrow 2.$ is exactly given by Proposition~\ref{Proposition: CC w-unital <=> EE unital}.

\noindent $2. \Rightarrow 3$.  If $\CC$ is w-unital, we can consider the free algebra $X = F(1)$ on one generator, and then the diagram
%\begin{equation}\label{coproduct}
$$
\xymatrix@1{
X & X+X \ar[r]^-{\cop{0}{1}} \ar[l]_-{\cop{1}{0}} & X }
$$
%\end{equation}
is a weak product in $\CC$ by Proposition~\ref{Proposition: characterization of w-unital cats with +}. It follows that there exists a morphism $\add:\map{X}{X+X}$ with the property that $\cop{1}{0} \cdot \add = 1 =  \cop{0}{1} \cdot \add$, and $X$ is a unital coalgebra.

\noindent $3. \Rightarrow 2.$ 
If $X = F(1)$ is a unital coalgebra, then \emph{every} free algebra $A$ is a unital coalgebra (since $A$ is a copower $A = \coprod_{i\in I} X$ for some set $I$).
Consider then two morphisms $b :\map{A}{B}$ and $c:\map{A}{C}$ in $\mathbb C$. Then the unital coalgebra structure on $A$ allows one to define a morphism 
$$\xymatrix{ A \ar[r]^-{\add} &A+A \ar[r]^-{(b+c)} & B+C}$$
such that $\cop{1}{0}\cdot (b+c)\cdot \add=b$ and $\cop{0}{1}\cdot(b+c)\cdot \add=c$. This proves that the span 
$$
\xymatrix@1{
B & B+C \ar[r]^-{\cop{0}{1}} \ar[l]_-{\cop{1}{0}} & C }
$$
is a weak product,  and $\CC$ is w-unital.

\noindent $3. \Leftrightarrow 4$.  Let $\mathbb T$ be an algebraic theory of $\EE$  (in the sense of Lawvere \cite{Law}), namely a small category with finite products with the property that $\EE$ is the category of product preserving functors from $\CC$ to the category of sets. Via the duality between varieties of algebras and algebraic theories, the existence of a morphism $\add:\map{X}{X+X}$ making the diagram
$$
\xymatrix@1@C=30pt{
X \ar@{=}[dr] & X + X \ar[l]_-{\cop{1}{0}} \ar[r]^-{\cop{0}{1}} & X \ar@{=}[dl] \\
& X \ar@{.>}[u]_-{\add} & }
$$
commute in $\CC$ corresponds exactly to the existence of a morphism \\$+:\map{T \times T}{T}$ making the diagram
$$
\xymatrix@1@C=30pt{
T \ar@{=}[dr] \ar[r]^-{\langle 1,0 \rangle} & T \times T \ar@{.>}[d]^+ & \ar[l]_-{\langle 0,1 \rangle} T  \ar@{=}[dl] \\
& T & }
$$
commute in the algebraic theory $\mathbb T$ (here $T$ is the ``generic object'' of $ \mathbb T$). 
\hfill $\Box$\vspace{5pt}

The proof above of the equivalence between conditions 3. and 4. follows the technique used in \cite{CP} in the case of Mal'tsev varieties \cite{DCMC}.

%%%%%%%%%%%%%%%%%%%%%%%%%%%%%%%%%%%%%%%  ÊXAMPLE: JONSSON-TARSKI QUASIVARIETY  %%%%%%%%%%%%%%%%%%%%%%%%%%%%%%%%%%%%%%%%
\begin{remark}
The arguments used in Theorem \ref{Theorem: Jonsson-Tarski and (w-)unital} can be adapted also to the context of \emph{quasivarieties} \cite{AP} (we never used the fact that $\EE$ is an exact category, but only the fact that $\EE$ is regular).
An interesting example of a unital quasivariety is provided by the category of \emph{torsion-free monoids}, which are the monoids $M$ satisfying the implications $$a^n = b^n \Rightarrow a =b, \qquad \forall a, b \in M, \quad \forall n \in \mathbb{N}^*.$$
\end{remark}

%%%%%%%%%%%%%%%%%%%%%%%%%%%%%%%%%%%%%%%%%%%%%%%%%%%%%%%%%%%%%%%%%%%%%%%%%%%%%%%%%%%%%%%%%%%%%%%%%%%%%%%%%%%%%%%%%%%%%%%
%%%%%%%%%%%%%%%%%%%%%%%%%%%%%%%%%%%%%%%  Subtractive categories  %%%%%%%%%%%%%%%%%%%%%%%%%%%%%%%%%%%%%%%%%%%%%%%%%%%%%%
%%%%%%%%%%%%%%%%%%%%%%%%%%%%%%%%%%%%%%%%%%%%%%%%%%%%%%%%%%%%%%%%%%%%%%%%%%%%%%%%%%%%%%%%%%%%%%%%%%%%%%%%%%%%%%%%%%%%%%%
\section{Subtractive categories}

In this section we shall characterize the projective covers of regular subtractive categories \cite{SC}.
In the presence of finite coproducts, this will lead to a simple proof of the characterization of algebraic pointed subtractive (quasi)varieties, in the sense of Ursini \cite{OSV}.

%%%%%%%%%%%%%%%%%%%%%%%%%%%%%%%%%%%%%%%  RIGHT/LEFT PUNCTUAL REFLEXIVE RELATION  %%%%%%%%%%%%%%%%%%%%%%%%%%%%%%%%%%%%%%
Given a pointed category $\EE$, a span $\Span{X}{f}{Z}{g}{Y}$ in $\EE$ is called
\emph{right (resp. left) punctual} if there exists a morphism $t:\map{Y}{Z}$ (resp. $s:\map{X}{Z}$) such that $g\cdot t=1$
and $f\cdot t=0$ (resp. $f\cdot s=1$ and $g\cdot s=0$). Recall that a reflexive graph is a span
$\punctualspan{X}{d}{e}{G}{c}{e}{X}$ such that $d\cdot e=c\cdot e=1$; a reflexive relation is a reflexive graph such
that $d$ and $c$ are jointly monomorphic. In the text, we shall omit representing the reflexivity morphism $e$ for reflexive graphs
and relations.

%%%%%%%%%%%%%%%%%%%%%%%%%%%%%%%%%%%%%%%  DEFINITION: SUBTRACTIVE CATEGORY  %%%%%%%%%%%%%%%%%%%%%%%%%%%%%%%%%%%%%%%%%%%%
\begin{definition}\emph{(\cite{SC})
\label{Subtractive}
A category $\EE$ is called \emph{subtractive} when it is pointed, has finite limits and every right punctual reflexive
relation is also left punctual.}
\end{definition}

Equivalently, a subtractive category could be defined by demanding that every left punctual reflexive relation
is also right punctual. We now introduce a ``weak version'' of this definition.

%%%%%%%%%%%%%%%%%%%%%%%%%%%%%%%%%%%%%%%  DEFINITION: w-SUBTRACTIVE CATEGORY  %%%%%%%%%%%%%%%%%%%%%%%%%%%%%%%%%%%%%%%%%%
\begin{definition}
\label{Definition: w-subtractive cat}
\emph{A category $\CC$ is called \emph{w-subtractive} when it is pointed, has weak finite limits and every right punctual
reflexive graph is also left punctual.}
\end{definition}

%%%%%%%%%%%%%%%%%%%%%%%%%%%%%%%%%%%%%%%  PROPOSITION: CC w-SUBTRACTIVE <=> EE SUBTRACTIVE  %%%%%%%%%%%%%%%%%%%%%%%%%%%%
\begin{proposition}
\label{Proposition: CC w-sutractive <=> EE subtractive}
Let $\CC$ be a pointed projective cover of a regular category $\EE$. Then $\EE$ is subtractive if and only if $\CC$ is
w-subtractive.
\end{proposition}
\proof Suppose that $\EE$ is subtractive. Since $\CC$ is a projective cover of a regular category $\EE$, then $\CC$ has weak finite limits (Remark~\ref{Remark: Proposition 4 in REC}). Now, let
$\rightpunctualspan{\tilde X}{\tilde d}{}{\tilde G}{\tilde c}{\tilde t}{\tilde X}$ be a right punctual reflexive graph
in $\CC$. Then, the (regular epimorphism, monomorphism)-factorization $(\tilde{d},\tilde{c})=(r_1,r_2)\cdot p$ in
$\EE$, gives a right punctual reflexive relation $\prd{r_1}{r_2}:\mono{R}{\tilde{X}\times \tilde{X}}$ in $\EE$. So, it is
left punctual, i.e. there exists a morphism $s:\map{\tilde{X}}{R}$ such that $r_1\cdot s=1$ and $r_2\cdot
s=0$. By the projectivity of $\tilde{X}$, there exists a morphism
$\tilde{s}:\map{\tilde{X}}{\tilde{G}}$ such that $p\cdot \tilde{s}=s$. This morphism $\tilde{s}$ makes the original
reflexive graph also left punctual.

\noindent Conversely, $\EE$ is pointed, since $\CC$ is (see Section~\ref{Projective covers}). Let
$\rightpunctualspan{X}{r_1}{}{R}{r_2}{t}{X}$ be a right punctual reflexive relation in $\EE$. Given  $\CC$-covers
$x:\regepi{\tilde X}{X}$ and $\alpha:\regepi{\tilde{R}}{x^{-1}(R)}$,
$$
  \xymatrix@1@R=25pt{\tilde R \ar@{>>}[r]^-{\alpha} \ar[dr]_-{\langle \tilde{r}_1, \tilde{r}_2 \rangle}
  & x^{-1}(R) \ar@{>>}[r]^-{u} \ar[d]|(.4){\langle r_1',r_2' \rangle\;\;} \ar@{}[dr]|(.25){\lrcorner} & R \ar[d]^-{\langle r_1,r_2 \rangle} \\
  & \tilde X\times \tilde X \ar@{>>}[r]_-{x\times x} & X\times X }
$$
the inverse image $\rightpunctualspan{\tilde{X}}{r_1'}{}{x^{-1}(R)}{r_2'}{t'}{\tilde{X}}$ is necessarily a reflexive
relation in $\EE$ which is also right punctual for $t'=(\prd{0}{1},t\cdot x)$. Using the projectivity of $\tilde{X}$ and the fact that
$\alpha$ is a regular epimorphism, there exists a morphism $\tilde{t}:\map{\tilde{X}}{\tilde{R}}$ such that
$\rightpunctualspan{\tilde{X}}{\tilde{r}_1}{}{\tilde{R}}{\tilde{r}_2}{\tilde{t}}{\tilde{X}}$ is a right punctual
reflexive graph in $\CC$. Therefore, it is also left punctual, i.e. there exists a morphism
$\tilde{s}:\map{\tilde{X}}{\tilde{R}}$ such that $\tilde{r}_1\cdot \tilde{s}=1$ and $\tilde{r}_2\cdot
\tilde{s}=0$. Finally, the fact that $x$ is a strong (=regular) epimorphism, gives a unique factorization $s$ in
$$
  \xymatrix@1@C=30pt{ \tilde{X} \ar@{>>}[r]^-x \ar[d]_-{u\cdot \alpha\cdot \tilde{s}} & X \ar[d]^-{\langle 1,0 \rangle} \ar@{.>}[dl]_-{s} \\
  R \ar@{ >->}[r]_-{\langle r_1,r_2 \rangle} & X\times X}
$$
which makes the original reflexive relation $R$ in $\EE$ also left punctual.
$\Box$\vspace{5pt}

%%%%%%%%%%%%%%%%%%%%%%%%%%%%%%%%%%%%%%%  FINITE +  %%%%%%%%%%%%%%%%%%%%%%%%%%%%%%%%%%%%%%%%%%%%%%%%%%%%%%%%%%%%%%%%%%%%
In the presence of binary coproducts, w-subtractive categories can be characterized through a special kind of right punctual
reflexive graph as follows:

%%%%%%%%%%%%%%%%%%%%%%%%%%%%%%%%%%%%%%%  PROPOSITION: CHARACTERIZATION OF w-SUBTRACTIVE CATEGORIES WITH +  %%%%%%%%%%%%
\begin{proposition}
\label{Proposition: characterization of w-subtractive cats with +}
A pointed weakly lex category with binary coproducts $\CC$ is w-sub\-trac\-tive if and only, for any object $X$ in $\CC$,
the right punctual reflexive graph $\rightpunctualspan{X}{\cop{1}{0}}{}{X+X}{\cop{1}{1}}{i_2}{X}$ is left punctual.
\end{proposition}
\proof The statement is obvious when $\CC$ is w-subtractive. Conversely, consider a right punctual reflexive graph
$\rightpunctualspan{X}{d}{}{G}{c}{t}{X}$. By assumption, there exists a morphism $s:\map{X}{X+X}$ such that $\cop{1}{0}\cdot s=1$
and $\cop{1}{1}\cdot s=0$. Consequently, there exists a morphism $\cop{e}{t}\cdot s:\map{X}{G}$, such that $d\cdot \cop{e}{t}\cdot s=1$ and
$c\cdot \cop{e}{t}\cdot s=0$, proving that the reflexive graph is also left punctual.
$\Box$\vspace{5pt}

Recall that a subtractive variety $\EE$, in the sense of Ursini \cite{OSV}, is such that its theory contains a 
constant 0 and a binary operation $s$, called \emph{subtraction}, such that $s(x,0)=x$ and $s(x,x)=0$.
When in the theory of a variety there is a \emph{unique} constant, then $\EE$ is a pointed variety. In this case it is known that $\EE$ is a subtractive variety if and only if
$\EE$ is a subtractive category (Theorem 2 of \cite{SC}).

The morphism $s:\map{X}{X+X}$ arising in the proof of Proposition \ref{Proposition: characterization of w-subtractive cats with +} then equips $X$ with a \emph{subtractive coalgebra} structure, i.e. such that the equalities $\cop{1}{0}\cdot s=1$ and $\cop{1}{1}\cdot s=0$ hold.

%%%%%%%%%%%%%%%%%%%%%%%%%%%%%%%%%%%%%%%  COROLLARY: CHARACTERIZATION OF w-SUBTRACTIVE WITH SUBTRACTIVE COALGEBRAS  %%%%
\begin{corollary}
\label{Corollary: characterization of w-subtractive cats with coalgs}
A pointed weakly lex category with binary coproducts $\CC$ is w-sub\-trac\-tive if and only every object $X$ in $\CC$ is a subtractive coalgebra.
\end{corollary}

%%%%%%%%%%%%%%%%%%%%%%%%%%%%%%%%%%%%%%%  SUBTRACTIVE VARIETIES  %%%%%%%%%%%%%%%%%%%%%%%%%%%%%%%%%%%%%%%%%%%%%%%%%%%%%%%

We are then ready to prove the main result in this section (see the proof of Theorem~\ref{Theorem: Jonsson-Tarski and (w-)unital} for the omitted details which we avoid repeating):

%%%%%%%%%%%%%%%%%%%%%%%%%%%%%%%%%%%%%%%  THEOREM: SUBTRACTIVE ISSUES  %%%%%%%%%%%%%%%%%%%%%%%%%%%%%%%%%%%%%%%%%%%%%%%%%
\begin{theorem}
\label{Theorem: Subtractive issues}
Let $\EE$ be an algebraic variety and $\CC$ its full subcategory of free algebras. Then following statements are
equivalent:
\begin{itemize}
  \item[1.] $\EE$ is a subtractive category;
  \item[2.] $\CC$ is a w-subtractive category;
  \item[3.] the free algebra $X=F(1)$ on one generator is a subtractive coalgebra;
  \item[4.] $\EE$ is a pointed subtractive variety.
\end{itemize}
\end{theorem}
\proof The conditions $1.$  and $2.$ are equivalent by Proposition~\ref{Proposition: CC w-sutractive <=> EE
subtractive}.  \\ The conditions $2.$ and $3.$ are equivalent by Corollary~\ref{Corollary: characterization of w-subtractive cats with coalgs}.

\noindent $3. \Leftrightarrow 4$. By duality, the span $$\punctualspan{X}{\cop{1}{0}}{s}{X+X}{\cop{1}{1}}{i_2}{X}$$ is punctual in $\CC$ if and only if the corresponding span
$$\longpunctualspan{T}{s}{\langle 1,0 \rangle}{T\times T}{\pi_2}{\langle 1,1 \rangle}{T}$$
is punctual in the algebraic theory $\mathbb T$ of $\EE$. \hfill 
$\Box$\vspace{5pt}

%%%%%%%%%%%%%%%%%%%%%%%%%%%%%%%%%%%%%%%%%%%%%%%%%%%%%%%%%%%%%%%%%%%%%%%%%%%%%%%%%%%%%%%%%%%%%%%%%%%%%%%%%%%%%%%%%%%%%%%
%%%%%%%%%%%%%%%%%%%%%%%%%%%%%%%%%%%%%%%  Strongly unital categories  %%%%%%%%%%%%%%%%%%%%%%%%%%%%%%%%%%%%%%%%%%%%%%%%%%
%%%%%%%%%%%%%%%%%%%%%%%%%%%%%%%%%%%%%%%%%%%%%%%%%%%%%%%%%%%%%%%%%%%%%%%%%%%%%%%%%%%%%%%%%%%%%%%%%%%%%%%%%%%%%%%%%%%%%%%
\section{Strongly unital categories}
In this last section, we are going to examine the characterization for projective covers of regular strongly unital categories. The link with unital and subtractive categories is given by the categorical equation (see Proposition 3 of \cite{SC}):
\begin{center}
	unital + subtractive = strongly unital.
\end{center}

%%%%%%%%%%%%%%%%%%%%%%%%%%%%%%%%%%%%%%%  DEFINITION: STRONGLY UNITAL CATEGORY  %%%%%%%%%%%%%%%%%%%%%%%%%%%%%%%%%%%%%%%%
\begin{definition}\emph{(\cite{MCFPO})
\label{Strongly unital}
A category $\EE$ is called \emph{strongly unital} when it is pointed, has finite limits and for every object $X$ in
$\EE$, the pair of morphisms\\ $\xymatrix@1@C=30pt{X \ar[r]^-{\langle 1,1 \rangle} & X\times X &  X \ar[l]_-{\langle 0,1 \rangle}}$ is jointly strongly
epimorphic.}
\end{definition}

%%%%%%%%%%%%%%%%%%%%%%%%%%%%%%%%%%%%%%%  SPLIT RIGHT PUNCTUAL SPAN  %%%%%%%%%%%%%%%%%%%%%%%%%%%%%%%%%%%%%%%%%%%%%%%%%%%
Let $\EE$ be a pointed category. A span $\Span{X}{f}{Z}{g}{Y}$ in $\EE$ is called \emph{split right punctual} if there
exist morphisms $s:\map{X}{Z}$ and $t:\map{Y}{Z}$ such that $f\cdot s=1, g\cdot t=1$ and $f\cdot t=0$, i.e. it is right
punctual and $f$ is split by $s$. 
For example, for any object $X$ in $\EE$, we have a split right punctual span
$\longpunctualspan{X}{\pi_1}{\langle 1,1 \rangle}{X\times X}{\pi_2}{\langle 0,1 \rangle}{X}$. It is actually possible to characterize strongly unital categories through such spans (see Proposition 1.8.14 \cite{BB}). This fact led us to introduce the following definition:

%%%%%%%%%%%%%%%%%%%%%%%%%%%%%%%%%%%%%%%  DEFINITION: STRONGLY w-UNITAL CATEGORY  %%%%%%%%%%%%%%%%%%%%%%%%%%%%%%%%%%%%%%
\begin{definition}
\label{Definition: strongly w-unital cat}
\emph{A category $\CC$ is called \emph{w-strongly unital} when it is pointed, has weak finite limits and, for every split
right punctual span $\punctualspan{X}{f}{s}{Z}{g}{t}{Y}$, the span $\Span{X}{f}{Z}{g}{Y}$ is a weak product.}
\end{definition}

The ``equation'' relating the notions of (strongly) unital and of subtractive categories still holds in the weakened context:

%%%%%%%%%%%%%%%%%%%%%%%%%%%%%%%%%%%%%%%  PROPOSITION: STRONGLY w-UNITAL = w-UNITAL + w-SUBTRACTIVE  %%%%%%%%%%%%%%%%%%%
\begin{proposition}
\label{Proposition: w-strongly unital = w-unital + w-subtractive}
A pointed weakly lex category $\CC$ is w-strongly unital if and only if it is w-unital and w-subtractive.
\end{proposition}
\proof Let $\CC$ be a w-strongly unital category. Since any punctual span is a split right punctual span, it then determines a
weak product, and the category $\CC$ is w-unital. To prove that $\CC$ is w-subtractive, consider a right punctual reflexive graph $\rightpunctualspan{X}{d}{}{G}{c}{t}{X}$. It is a split right
punctual span, thus a weak product. Then, there exists a morphism $s:\map{X}{G}$ such that $d\cdot s=1$ and $c\cdot s=0$,
which makes the reflexive graph also left punctual. 

\noindent Conversely, suppose that $\CC$ is w-unital and w-subtractive. Let
$\punctualspan{X}{f}{s}{Z}{g}{t}{Y}$ be a split right punctual span and consider arbitrary morphisms $x:\map{A}{X}$ and
$y:\map{A}{Y}$. Any weak pullback of $\prd{f}{g}$ along $f\times g$
$$
\xymatrix@1{
  G \ar[r]^-h \ar[d]_-{\langle d,c \rangle} & Z \ar[d]^-{\langle f,g \rangle} \\ Z\times Z \ar[r]_-{f\times g} & X\times Y}
$$
produces a right punctual \emph{reflexive} graph $\rightpunctualspan{Z}{d}{}{G}{c}{\tau}{Z}$, where the reflexive morphism is
$e=(\prd{1}{1},1)$ and $\tau=(\prd{0}{t},t)$. Then, it must also be left punctual and, consequently, it is a weak product. This
gives a morphism $w:\map{A}{G}$ such that $d\cdot w=s\cdot x$ and $c\cdot w=t\cdot y$. Then, there exists a morphism
$h\cdot w:\map{A}{Z}$ such that $f\cdot h\cdot w=x$ and $g\cdot h\cdot w=y$, proving that $(f,g)$ is a weak product.
\hfill $\Box$\vspace{5pt}

%%%%%%%%%%%%%%%%%%%%%%%%%%%%%%%%%%%%%%%  PROPOSITION: CC STRONGLY w-UNITAL <=> EE STRONGLY UNITAL  %%%%%%%%%%%%%%%%%%%%
\begin{corollary}
\label{Proposition: CC strongly w-unital <=> EE strongly unital}
Let $\CC$ be a pointed projective cover of a regular category $\EE$. Then $\EE$ is strongly unital if and only if
$\CC$ is w-strongly unital.
\end{corollary}
\proof The proof follows from Proposition~\ref{Proposition: w-strongly unital = w-unital + w-subtractive} and Proposition $3$ in \cite{SC}. 
$\Box$\vspace{5pt}

%%%%%%%%%%%%%%%%%%%%%%%%%%%%%%%%%%%%%%%  FINITE +  %%%%%%%%%%%%%%%%%%%%%%%%%%%%%%%%%%%%%%%%%%%%%%%%%%%%%%%%%%%%%%%%%%%%
In a context that admits binary coproducts, (w-)strongly unital categories can also be characterized through a special kind of
punctual span. For strongly unital categories, the characterization is similar to that of Proposition 1.2.18 in \cite{BB}. We omit the proof of the following result, which is similar to the one of Proposition~\ref{Proposition: characterization of w-unital cats with +}.

%%%%%%%%%%%%%%%%%%%%%%%%%%%%%%%%%%%%%%%  PROPOSITION: CHARACTERIZATION OF STRONGLY w-UNITAL CATEGORIES WITH +  %%%%%%%%
\begin{proposition}
\label{Proposition: characterization of strongly w-unital cats with +}
A pointed weakly lex category with binary coproducts $\CC$ is w-strongly unital if and only the span
$\longspan{X}{\cop{1}{0}}{X+X}{\cop{1}{1}}{X}$ is a weak product, for every object $X$ in $\CC$.
\end{proposition}
%\proof The statement is obvious when $\CC$ is strongly w-unital, since the span
%$\longpunctualspan{X}{\cop{1}{0}}{i_1}{X+X}{\cop{1}{1}}{i_2}{X}$ is a split right punctual span. Conversely, consider a split
%right punctual span $\punctualspan{X}{f}{s}{Z}{g}{t}{Y}$ and arbitrary morphisms $x:\map{A}{X}, y:\map{A}{Y}$. By
%assumption, there exists a morphism $w:\map{A}{Z+Z}$ such that $\cop{1}{0}\cdot w=s\cdot x$ and $\cop{1}{1}\cdot w=t\cdot y$. So,
%there exists a morphism $\cop{1}{t\cdot g}\cdot w:\map{A}{Z}$ such that $f\cdot \cop{1}{t\cdot g}\cdot w=x$ \and $g\cdot \cop{1}{t\cdot g}\cdot w=y$, %proving that $(f,g)$ is a weak product.
%$\Box$\vspace{5pt}

%%%%%%%%%%%%%%%%%%%%%%%%%%%%%%%%%%%%%%%  POINTED STRONGLY UNITAL VARIETIES  %%%%%%%%%%%%%%%%%%%%%%%%%%%%%%%%%%%%%%%%%%%
Recall that a pointed variety is strongly unital if and only if its corresponding theory contains a unique constant 0 and a
ternary operation $p$ such that $p(x,0,0)=x$ and $p(x,x,y)=y$ (Theorem 1.8.16 of \cite{BB}). This result is illustrated in
the next theorem, that uses the w-strongly unital property. An object $X$ equipped with a morphism $p:\map{X}{X+X+X}$ such that $\left<\begin{smallmatrix} 1\\ 0\\0 \end{smallmatrix}\right>\cdot p=1$ and $(\cop{1}{1}+1)\cdot p=i_2$ is called a \emph{strongly unital coalgebra}.

%%%%%%%%%%%%%%%%%%%%%%%%%%%%%%%%%%%%%%%  THEOREM: STRONLGY (w-)UNITAL  %%%%%%%%%%%%%%%%%%%%%%%%%%%%%%%%%%%%%%
\begin{theorem}
\label{Theorem: strongly (w-)unital}
Let $\EE$ be an algebraic variety and $\CC$ its full subcategory of free algebras. Then following statements are
equivalent:
\begin{itemize}
  \item[1.] $\EE$ is a strongly unital category;
  \item[2.] $\CC$ is a strongly w-unital category;
  \item[3.] the free algebra $X=F(1)$ on one generator is a strongly unital coalgebra;
  \item[4.] $\EE$ is a pointed strongly unital variety.
\end{itemize}
\end{theorem}
\proof $1. \Leftrightarrow 2.$ is given by Corollary~\ref{Proposition: CC strongly w-unital <=> EE
strongly unital}. 

\noindent $2. \Rightarrow 3$. If $\CC$ is w-strongly unital, then for the free algebra $X=F(1)$, the split
right punctual span $\longpunctualspan{X}{\left<\begin{smallmatrix} 1\\ 0\\0 \end{smallmatrix}\right>}{i_1}{X+X+X}{\cop{1}{1}+1}{i_{2,3}}{X+X}$ gives a weak product. Therefore, there exists a morphism $p:\map{X}{X+X+X}$ such that $\left<\begin{smallmatrix} 1\\ 0\\0 \end{smallmatrix}\right>\cdot p=1$ and $(\cop{1}{1}+1)\cdot p=i_2$; $X$ is a strongly unital coalgebra.

\noindent $3. \Rightarrow 2.$ By using the same argument as in the implication $3. \Rightarrow 2.$ of Theorem \ref{Theorem: Jonsson-Tarski and (w-)unital} we know that any free algebra $A$ has a strongly unital coalgebra structure since the free algebra $X= F(1)$ has such a structure. Consider then free algebra morphisms $b, c:\map{A}{Y}$. There exists a morphism 

$$\xymatrix@C=40pt{A \ar[r]^-p & A+A+A \ar[r]^-{b+\cop{b}{c}} & Y+Y
}$$
such that $\cop{1}{0}\cdot \left(b+\cop{b}{c}\right)\cdot p=b$ and $\cop{1}{1}\cdot \left(b+\cop{b}{c}\right)\cdot p=c$. This proves that the span 
$\longspan{Y}{\cop{1}{0}}{Y+Y}{\cop{1}{1}}{Y}$
 is a weak product for any $Y \in \CC$, and $\CC$ is w-strongly unital by Proposition~\ref{Proposition: characterization of strongly w-unital cats with +}.

\noindent $3. \Leftrightarrow 4$.  Via the duality between varieties of algebras and algebraic theories, the existence of a morphism $p:\map{X}{X+X+X}$ making the diagram
$$
\xymatrix@1@C=35pt{
X \ar@{=}[dr] & X + X + X \ar[l]_-{\left<\begin{smallmatrix} 1\\ 0\\0 \end{smallmatrix}\right>} \ar[r]^-{\cop{1}{1}+1} & X + X \\
& X \ar@{.>}[u]_p \ar[ur]_-{i_2} & }
$$
 commutative
in $\CC$ (where $X=F(1)$) corresponds to the existence of a morphism $p:\map{T\times T\times T}{T}$ making the diagram
$$
\xymatrix@1@C=40pt{
T \ar@{=}[dr] \ar[r]^-{\langle 1,0,0 \rangle} & T\times T\times T \ar@{.>}[d]^p & \ar[l]_-{\langle 1,1 \rangle\times 1} T  \ar[dl]^-{\pi_2} \\
& T & }
$$
commutative in the algebraic theory $\mathbb T$ of $\EE$.

\hfill $\Box$
%%%%%%%%%%%%%%%%%%%%%%%%%%%%%%%%%%%%%%%%%%%%%%%%%%%%%%%%%%%%%%%%%%%%%%%%%%%%%%%%%%%%%%%%%%%%%%%%%%%%%%%%%%%%%%%%%%%%%%%
%%%%%%%%%%%%%%%%%%%%%%%%%%%%%%%%%%%%%%%  REFERENCES  %%%%%%%%%%%%%%%%%%%%%%%%%%%%%%%%%%%%%%%%%%%%%%%%%%%%%%%%%%%%%%%%%%
%%%%%%%%%%%%%%%%%%%%%%%%%%%%%%%%%%%%%%%%%%%%%%%%%%%%%%%%%%%%%%%%%%%%%%%%%%%%%%%%%%%%%%%%%%%%%%%%%%%%%%%%%%%%%%%%%%%%%%%


\begin{thebibliography}{NEKEAC}
\bibitem{AP} J. Ad\'amek and H. Porst, \emph{Algebraic theories of quasivarieties}, J. Algebra, \textbf{208} (1998) 379-398.
\bibitem{Adelman} M. Adelman, \emph{Abelian categories over additive ones}, J. Pure Appl. Algebra \textbf{3} (1973) 103-117.
\bibitem{EC} M. Barr, Exact Categories, in: \emph{Lecture Notes in Math.} 236, Springer (1971) 1-120.
\bibitem{BB} F. Borceux and D. Bourn, \emph{Mal'cev, protomodular, homological and semi-abelian categories}, Kluwer, 2004.
\bibitem{MCFPO} D. Bourn, \emph{Mal'cev categories and fibration of pointed objects}, Appl. Categ. Structures \textbf{4} (1996) 307-327.
%\bibitem{ST} M. Bunge, A. Carboni, \emph{The symmetric topos}, J. Pure Appl. Algebra \textbf{105} (1995) 233-249.
\bibitem{DCMC} A. Carboni, J. Lambek, and M.C. Pedicchio, \emph{Diagram chasing in Malcev categories}, J. Pure Appl. Algebra \textbf{69} (1991) 271-284.
\bibitem{FECLEO} A. Carboni and R. Celia Magno, \emph{The free exact category on a left exact one}, J. Austr. Math. Soc. Ser., A 33 (1982) 295-301.
\bibitem{CP} A. Carboni and M.C. Pedicchio, \emph{A new proof of the Mal'cev theorem}, Rend. Circolo Mat. Palermo, S. II, Suppl. \textbf{64} (2000) 13-16.
\bibitem{LCCEC} A. Carboni and G. Rosolini, \emph{Locally cartesian closed exact completions}, J. Pure Appl. Algebra \textbf{154} (2000) 103-116.
\bibitem{REC} A. Carboni and E. Vitale, \emph{Regular and exact completions}, J. Pure Appl. Algebra \textbf{125} (1998) 79-116.
\bibitem{G} M. Gran, \emph{Semi-abelian exact completions}, Homology, Hom. Appl. \textbf{4} (2002) 175-189.
\bibitem{ECHC} M. Gran and E.M. Vitale, \emph{On the exact completion of the homotopy category}, Cah. Top. Geom. Diff. Cat´egoriques \textbf{XXXIX-4} (1998) 287-297.
\bibitem{RAC} P. Freyd, \emph{Representations in abelian categories}, in: Proc. Conference on Categorical Algebra, La Jolla,
1965 (Springer, Berlin, 1966).
\bibitem{SC} Z. Janelidze, \emph{Subtractive categories}, Appl. Categ. Structures \textbf{13} (2005) 343-350.
\bibitem{Zurab} Z. Janelidze,
\emph{Closedness properties of internal relations VI: Approximate operations},
Cah. Topol. et G\'eom. Diff\'er. Cat\'eg. 50 (2009)
298-319.
\bibitem{JT} B. J\'onsson and A. Tarski, \emph{Direct Decompositions of Finite Algebraic Systems}, Notre
Dame Mathematical Lectures, Notre Dame, Indiana (1947)
\bibitem{Law} F. W. Lawvere, \emph{Functorial Semantics of Algebraic Theories}, Ph.D. thesis, Columbia University (1963), Reprint in Theory and Applications of Categories \textbf{5} (2004) 1-121.
\bibitem{Menni} M. Menni, \emph{A characterization of the left exact categories whose exact completions are toposes}, J. Pure Appl. Algebra \textbf{177} (2003) 287-301.
%\bibitem{PV} M.C. Pedicchio and E.M. Vitale, \emph{On the abstract characterisation of quasi-varieties}, Algebra Universalis  \textbf{31} (2000) %269-278.
%\bibitem{Pitts} A. Pitts, \emph{The lex reflection of a category with finite products}, unpublished notes (1996).
\bibitem{CCEC} J. Rosicky, \emph{Cartesian closed exact completions}, J. Pure Appl. Algebra \textbf{142} (1999) 261-270.
\bibitem{ECRAC} J. Rosicky and E.M. Vitale, \emph{Exact completions and representations in abelian categories}, Homology, Hom. Appl. \textbf{3} (2001) 453-466.
\bibitem{OSV} A. Ursini, \emph{On subtractive varieties}, I. Algebra Universalis \textbf{31} (1994) 204-222.
\bibitem{Vit} E. Vitale, \emph{Left covering functors}, PhD thesis, Universit\'e catholique de Louvain (1994)
\end{thebibliography}
\end{document}